\input epsf
\documentclass{amsart}
\usepackage{amssymb,latexsym}

\numberwithin{equation}{section}
\newtheorem{theorem}{Theorem}[section]
\newtheorem{proposition}[theorem]{Proposition}
\newtheorem{corollary}[theorem]{Corollary}
\newtheorem{definition}[theorem]{Definition}
\newtheorem{conjecture}[theorem]{Conjecture}
\newtheorem{remark}[theorem]{Remark}
\newtheorem{lemma}[theorem]{Lemma}

\newtheorem{example}[theorem]{Example}

\def\ZZ{\mathbb{Z}}

\def\RR{\mathbb{R}}

\def\FF{\mathbb{F}}
\def\RR{\mathbb{R}}

\def\ii{\mathbf{i}}

\def\veps{\varepsilon}

\def\l{\ell}

\def\<{\langle}
\def\>{\rangle}

\newcommand{\beal}{\begin{eqnarray}\begin{array}{l} }
\newcommand{\bear}{\begin{eqnarray}\begin{array}{r} }
\newcommand{\beac}{\begin{eqnarray}\begin{array}{c} }
\newcommand{\bealn}{\begin{eqnarray*}\begin{array}{l} }
\newcommand{\bearn}{\begin{eqnarray*}\begin{array}{r} }
\newcommand{\beacn}{\begin{eqnarray*}\begin{array}{c} }
\newcommand{\eea}{\end{array}\end{eqnarray}}
\newcommand{\eean}{\end{array}\end{eqnarray*}}

\newcommand{\beq}{\begin{equation} }
\newcommand{\eeq}{\end{equation} }

\newcommand{\mat}[4]{\left(\begin{array}{cc}
#1 & #2 \\
#3 & #4 \\
\end{array}\right)}

\newcommand{\doublesubscript}[3]{
\displaystyle\mathop{\displaystyle #1_{#2}}_{#3}
}

\begin{document}

\title{Simply-laced Coxeter groups and groups generated by symplectic
transvections}

\author{Boris Shapiro}
\address{Department of Mathematics, University of Stockholm, S-10691, Sweden}
\email{shapiro@matematik.su.se}

\author{Michael Shapiro}
\address{Department of Mathematics, Royal Institute of Technology, S-10044,
Sweden}
\email{mshapiro@math.kth.se}

\author{Alek Vainshtein}
\address{Departments of Mathematics and Computer Science, University of Haifa,
Mount Carmel,
31905 Haifa, Israel}
\email{alek@mathsc11.haifa.ac.il}

\author{Andrei Zelevinsky}
\address{Department of Mathematics,
Northeastern University, Boston, Massachusetts 02115, USA}
\email{andrei@neu.edu}

\thanks{The research of Andrei Zelevinsky was supported in part 
by NSF grant \#DMS-9625511. 
The research of Boris Shapiro was supported in part
by NSF grants 
during his visit to Northeastern University
in the Fall 1997.}

\date{\today }


\maketitle

\makeatletter
\renewcommand{\@evenhead}{\tiny \thepage \hfill  B.~SHAPIRO, M.~SHAPIRO,
A.~VAINSHTEIN, A.~ZELEVINSKY \hfill}
\renewcommand{\@oddhead}{\tiny \hfill  SIMPLY-LACED COXETER GROUPS
 \hfill \thepage}
\makeatother

\section{Introduction}
\label{sec:intro}

The point of departure for this paper is the following result obtained
in \cite{ssv1, ssv2}.
Let $N_n^0$ denote the semi-algebraic set of all
unipotent upper-triangular $n \times n$ matrices $x$ with real entries
such that, for every $k = 1, \ldots, n-1$,
the minor of $x$ with rows $1, \ldots, k$ and columns
$n-k+1, \ldots, n$ is non-zero.
Then the number $\#_n$ of connected components of $N_n^0$ is given
as follows: $\#_2 = 2, \#_3 = 6, \#_4 = 20, \#_5 = 52$, and
$\#_n = 3 \cdot 2^{n-1}$ for $n \geq 6$.

An interesting feature of this answer is that every case which one can 
check by hand turns out to be exceptional.
But the method of the proof seems to be even more interesting
than the answer itself: it is shown that the connected components
of $N_n^0$ are in a bijection with the orbits of a certain group
$\Gamma_n$ that acts in a
vector space of dimension $n (n-1)/2$ over the two-element field $\FF_2$,
and is generated by symplectic transvections.
Such groups appeared earlier in singularity theory, see e.g.,
\cite{Janssen} and references therein.

The construction of $\Gamma_n$ given in \cite{ssv1, ssv2} uses the 
combinatorial machinery 
(developed in \cite{BFZ})
of pseudo-line arrangements associated with reduced expressions
in the symmetric group.
In this paper we present the following
far-reaching generalization of this construction.
Let $W$ be an arbitrary Coxeter group of simply-laced type
(possibly infinite but of finite rank). 
Let $u$ and $v$ be any two elements in $W$, and $\ii$
be a reduced word (of length $m = \l (u) + \l (v)$)
for the pair $(u,v)$ in the Coxeter group $W \times W$
(see Section~\ref{sec:definitions} for more details).
We associate to $\ii$ a
subgroup $\Gamma_\ii$ in $GL_{m} (\ZZ)$
generated by symplectic transvections.
We prove among other things that the subgroups corresponding
to different reduced words for the same pair $(u,v)$
are conjugate to each other inside $GL_{m} (\ZZ)$.
To recover the group $\Gamma_n$ from this general construction,
one needs several specializations and reductions:
take $W$ to be the symmetric group $S_n$; take $(u,v) = (w_0, e)$,
where $w_0$ is the longest permutation in $S_n$, and
$e$ is the identity permutation;
take $\ii$ to be the lexicographically minimal reduced word
$1, 2, 1, \ldots, n-1, n-2, \ldots, 1$ for $w_0$; and finally,
take the group $\Gamma_\ii (\FF_2)$ obtained from $\Gamma_\ii$
by reducing the linear transformations from $\ZZ$ to $\FF_2$.

We also generalize the enumeration result of \cite{ssv1, ssv2} by showing that,
under certain assumptions on $u$ and $v$, the number of 
$\Gamma_\ii (\FF_2)$-orbits in 
$\FF_2^{m}$ is equal to $3 \cdot 2^s$, where $s$ is the number of 
simple reflections
in $W$ that appear in a reduced decomposition for $u$ or $v$. 
We deduce this from a description of orbits in an even more general situation  
which sharpens the results in \cite{Janssen,ssv2} (see 
Section~\ref{sec:proofs-enumeration} 
below). 

Although the results and methods of this paper are purely algebraic and 
combinatorial,
our motivation for the study of the groups $\Gamma_\ii$ and their orbits 
comes from geometry.
In the case when $W$ is the (finite) Weyl group of a simply-laced root system,
we expect that the $\Gamma_\ii (\FF_2)$-orbits in $\FF_2^{m}$ enumerate 
connected components of
the real part of the reduced double Bruhat cell corresponding to $(u,v)$.
Double Bruhat cells were introduced and studied in \cite{FZ} as 
a natural framework for the study of total positivity in semisimple groups; 
as explained to us by N.~Reshetikhin, they also appear naturally in the 
study of symplectic
leaves in semisimple groups  (see \cite{HKKR}). 
Let us briefly recall their definition.

Let $G$ be a split simply connected semisimple algebraic group defined over 
$\RR$
with the Weyl group $W$; thus $W = {\rm Norm}_G (H)/H$,
where $H$ is an $\RR$-split maximal torus in $G$.
Let $B$ and $B_-$ be two (opposite) Borel subgroups in $G$ such that
$B \cap B_- = H$.
The \emph{double Bruhat cells}~$G^{u,v}$ are defined as the intersections
of ordinary Bruhat cells taken with respect to $B$ and $B_-$:
$$G^{u,v} = B u B  \cap B_- v B_- \ .$$
In view of the well-known Bruhat decomposition,
the group $G$ is the disjoint union of all $G^{u,v}$ for
$(u,v) \in W \times W$.

The term ``cell" might be misleading because the topology
of $G^{u,v}$ can be quite complicated.
The torus $H$ acts freely on $G^{u,v}$ by left (as well as right) translations,
and there is a natural section $L^{u,v}$ for this
action which we call the \emph{reduced double Bruhat cell}. 
These sections are introduced and studied in a forthcoming paper \cite{BZ99}
(for the definition see Section~\ref{sec:components-claims} below).

We seem to be very close to a proof of the fact that the connected 
components 
of the real part of $L^{u,v}$ are in a natural bijection with 
the $\Gamma_\ii (\FF_2)$-orbits in $\FF_2^{m}$; but some details are still 
missing. 
This question will be treated in a separate publication.  

The special case when $(u,v) = (e,w)$ for some element $w \in W$
is of particular geometric interest.
In this case, $L^{u,v}$ is biregularly
isomorphic to the so-called \emph{opposite Schubert cell} 
\begin{equation*}
C_w^0 := C_w \cap w_0 C_{w_0} \ ,
\end{equation*}
where $w_0$ is the longest element of $W$, and
$C_w = (B w B)/B  \subset G/B$ is the \emph{Schubert cell}
corresponding to $w$.
These opposite cells appeared in the literature in various
contexts, and were studied (in various degrees of generality)
in \cite{BFZ,BZ,rie,rie2,ssv1,ssv2}.
In particular, the variety $N_n^0$ which was the main object of study in
\cite{ssv1, ssv2} is naturally identified with the real part of
the opposite cell $C_{w_0}^0$ for $G = SL_n$.

By the informal ``complexification principle" of Arnold,
if the group $\Gamma_\ii (\FF_2)$ enumerates connected components
of the real part of $L^{u,v}$, the group $\Gamma_\ii$ itself
(which acts in $\ZZ^{m}$ rather than in $\FF_2^{m}$)
should provide information about topology of the complex variety $L^{u,v}$.
So far we did not find a totally satisfactory
``complexification" along these lines.

The paper is organized as follows.
Main definitions, notations and conventions are collected in 
Section~\ref{sec:definitions}.
Our main results are formulated  in Section~\ref{sec:results} and proved in 
the three next sections. 
We conclude by discussing in more detail the geometric connection outlined 
above.

\section{Definitions}
\label{sec:definitions}

\subsection{Simply-laced Coxeter groups}
\label{sec:Coxeter groups general}

Let $\Pi$ be an arbitrary finite graph without loops and multiple edges.
Throughout the paper, we use the following notation: write $i \in \Pi$ if 
$i$ is a vertex
of $\Pi$, and $\{i,j\} \in \Pi$ if the vertices $i$ and $j$ are adjacent in 
$\Pi$.
The (simply-laced) Coxeter group $W = W(\Pi)$ associated with $\Pi$ is 
generated by the
elements $s_i$ for $i \in \Pi$ subject to the relations
\begin{equation}
\label{eq:Coxeter}
s_i^2 = e; \quad s_i s_j = s_j s_i \, \, (\{i,j\} \notin \Pi);
\quad s_i s_j s_i = s_j s_i s_j \, \, (\{i,j\} \in \Pi) \ .
\end{equation}
A word $\ii = (i_1, \ldots, i_m)$ in the alphabet $\Pi$
is a \emph{reduced word} for $w \in W$ if $w = s_{i_1} \cdots s_{i_m}$,
and $m$ is the smallest length of such a factorization.
The length $m$ of any reduced word for $w$ is called the \emph{length} of $w$
and denoted by $m = \l (w)$.
Let $R(w)$ denote the set of all reduced words for $w$.

The ``double" group $W \times W$ is also a Coxeter group; it
corresponds to the graph $\tilde \Pi$ which is the
union of two disconnected copies of $\Pi$.
We identify the vertex set of $\tilde \Pi$ with $\{+1, -1\} \times \Pi$,
and write a vertex $(\pm 1, i) \in \tilde \Pi$ simply as $\pm i$.
For each $\pm i \in \tilde \Pi$, we set $\veps (\pm i) = \pm 1$
and $|\pm i| = i \in \Pi$.
Thus two vertices $i$ and $j$ of $\tilde \Pi$ are joined by an edge
if and only if $\veps (i) = \veps (j)$ and $\{|i|, |j|\} \in \Pi$.
In this notation, a reduced word for a pair $(u,v) \in W \times W$
is an arbitrary shuffle of a reduced word for $u$ written in the alphabet
$-\Pi$ and a reduced word for $v$ written in the alphabet $\Pi$.

In view of the defining relations (\ref{eq:Coxeter}), the set
of reduced words $R(u,v)$ is equipped with the following operations:
\begin{itemize}
\item \emph{2-move}. Interchange two consecutive entries $i_{k-1},i_k$
in a reduced word $\ii = (i_1, \ldots, i_m)$ provided
$\{i_{k-1},i_k\} \notin \tilde \Pi$.
\item \emph{3-move}. Replace three consecutive entries $i_{k-2},i_{k-1},i_k$
in $\ii$ by $i_{k-1},i_{k-2},i_{k-1}$ if $i_k = i_{k-2}$ and
$\{i_{k-1},i_k\} \in \tilde \Pi$.
\end{itemize}
In each case, we will refer to the index $k \in [1,m]$ as the
\emph{position} of the corresponding move.
Using these operations, we make $R(u,v)$ the set of vertices of
a graph whose edges correspond to $2$- and $3$-moves.
It is a well known result due to Tits that this graph is \emph{connected}, 
i.e.,
any two reduced words in $R(u,v)$ can be obtained from each other by a
sequence of $2$- and $3$-moves.
We will say that a $2$-move interchanging the entries $i_{k-1}$ and $i_k$
is \emph{trivial} if $i_k \neq -i_{k-1}$;
the remaining $2$-moves and all $3$-moves will be referred to
as \emph{non-trivial}.

\subsection{Groups generated by symplectic transvections}
\label{sec:symplectic transvections general}

Let $\Sigma$ be a finite directed graph.
As before, we shall write $k \in \Sigma$ if $k$ is a vertex of $\Sigma$, 
and $\{k,l\} \in \Sigma$ if the vertices $k$ and $l$ are adjacent in 
the underlying graph obtained from $\Sigma$ by forgetting directions of edges. 
We also write $(k \to l) \in \Sigma$ if $k \to l$ is a directed edge of 
$\Sigma$.

Let $V = \ZZ^\Sigma$ be the lattice with a fixed $\ZZ$-basis 
$(e_k)_{k \in \Sigma}$
labeled by vertices of $\Sigma$. 
Let $\xi_k \in V^*$ denote the corresponding coordinate functions, i.e., 
every vector $v \in V$ can be written as 
$$v = \sum_{k \in \Sigma} \xi_k (v) e_k \ .$$
We define a skew-symmetric bilinear form $\Omega$ on $V$ by
\begin{equation}
\label{eq:Omega}
\Omega = \Omega_\Sigma = \sum_{(k \to l) \in \Sigma} \xi_k \wedge \xi_l \ .
\end{equation}

For each $k \in \Sigma$, we define the symplectic transvection
$\tau_k = \tau_{k,\Sigma}: V \to V$ by
\begin{equation}
\label{eq:transvections}
\tau_k (v) = v - \Omega(v, e_k) e_k \ .
\end{equation}
(The word ``symplectic" might be misleading since $\Omega$ is allowed to be 
degenerate;
still we prefer to keep this terminology from \cite{Janssen}.) 
In the coordinate form, we have $\xi_l (\tau_k (v)) = \xi_l (v)$ 
for $l \neq k$, and
\begin{equation}
\label{eq:transvections2}
\xi_k (\tau_k (v)) = \xi_k (v) - \sum_{(a \to k) \in \Sigma} \xi_a (v)+ 
\sum_{(k \to b) \in \Sigma} \xi_b (v) \ .
\end{equation}
For any subset $B$ of vertices of $\Sigma$, we denote by
$\Gamma_{\Sigma, B}$ the group of linear transformations of $V = \ZZ^\Sigma$
generated by the transvections $\tau_k$ for $k \in B$.

Note that all transformations from $\Gamma_{\Sigma, B}$ are represented by 
integer matrices
in the standard basis $e_k$. 
Let $\Gamma_{\Sigma, B}(\FF_2)$ denote the group of linear transformations 
of the 
$\FF_2$-vector space $V(\FF_2) = \FF_2^\Sigma$ obtained from 
$\Gamma_{\Sigma, B}$
by reduction modulo $2$ (recall that $\FF_2$ is the $2$-element field).

\section{Main results}
\label{sec:results}

\subsection{The graph $\Sigma(\ii)$}
\label{sec:main construction}
We now present our main combinatorial construction that brings together 
simply-laced Coxeter groups and groups generated by symplectic transvections. 
Let $W = W(\Pi)$ be the simply-laced Coxeter group associated to a graph $\Pi$
(see Section~\ref{sec:Coxeter groups general}).
Fix a pair $(u,v) \in W \times W$, and let $m = \l (u) + \l(v)$. 
Let $\ii = (i_1, \ldots, i_m) \in R(u,v)$ be any reduced word for $(u,v)$. 
We shall construct a directed graph $\Sigma (\ii)$ and a subset 
$B(\ii)$ of its vertices, thus giving rise to a group 
$\Gamma_{\Sigma (\ii), B(\ii)}$ generated by symplectic transvections. 

First of all, the set of vertices of $\Sigma (\ii)$ is just the set 
$[1,m] = \{1, 2, \ldots, m\}$.
For $l \in [1,m]$, we denote by $l^- = l^-_\ii$ the maximal index
$k$ such that $1 \leq k < l$ and $|i_k| = |i_l|$; if $|i_k| \neq |i_l|$
for $1 \leq k < l$ then we set $l^- = 0$. 
We define $B(\ii) \subset [1,m]$ as the subset of indices
$l \in [2,m]$ such that $l^- > 0$. 
The indices $l \in B(\ii)$ will be called \emph{$\ii$-bounded}.

It remains to define the edges of $\Sigma(\ii)$. 

\begin{definition}
\label{def:edges}
{\rm A pair $\{k,l\} \subset [1,m]$ with $k<l$ is an edge of $\Sigma(\ii)$ if 
it satisfies one of the following three conditions:

(i) $k=l^-$;

(ii) $k^- < l^- < k$, $\{|i_k|, |i_l|\} \in \Pi$, and
$\veps(i_{l^-})=\veps(i_{k})$;

(iii) $l^-< k^- < k$, $\{|i_k|, |i_l|\} \in \Pi$, and
$\veps(i_{k^-})=-\veps(i_{k})$.

\noindent The edges of type (i) are called \emph{horizontal}, 
and those of types (ii) and (iii) \emph{inclined}. 
A horizontal (resp. inclined) edge $\{k,l\}$ with $k < l$ is directed 
from $k$ to $l$ 
if and only if $\veps (i_{k}) = +1$ (resp. $\veps (i_{k}) = -1$).}
\end{definition}

We will give a few examples in the end of Section~\ref{sec:graphs}.

\subsection{Properties of graphs $\Sigma(\ii)$}
\label{sec:graphs}
We start with the following property of $\Sigma(\ii)$ and $B(\ii)$. 

\begin{proposition}
\label{pr:boundary vertex}
For any non-empty subset $S \subset B(\ii)$, there exists a 
vertex $a \in [1,m] \setminus S$
such that $\{a,b\} \in \Sigma (\ii)$ for a unique $b \in S$. 
\end{proposition}

For any edge $\{i,j\} \in \Pi$, let $\Sigma_{i,j}(\ii)$ denote the 
induced directed subgraph of $\Sigma({\ii})$ with vertices $k \in [1,m]$ 
such that $|i_k|=i$ or $|i_k|=j$. 
We shall use the following planar realization of $\Sigma_{i,j}(\ii)$ 
which we call the $(i,j)$-\emph{strip} of $\Sigma(\ii)$. 
Consider the infinite horizontal strip $\RR\times [-1,1] \subset\RR^2$, and 
identify each vertex $k \in \Sigma_{i,j}(\ii)$ with the point 
$A = A_k = (k,y)$, where $y = -1$ for $|i_k| = i$, and $y = 1$ for $|i_k| = j$.
We represent each (directed) edge $(k \to l)$ by a straight line segment 
from $A_k$ to $A_l$. (This justifies the terms ``horizontal" and 
``inclined"  edges in Definition~\ref{def:edges}.)

Note that every edge of $\Sigma({\ii})$ belongs to some $(i,j)$-strip,  
so we can think of $\Sigma({\ii})$ as the union of all its strips 
glued together along horizontal lines.

\begin{theorem}
\label{th:strip}

(a) The $(i,j)$-strip of $\Sigma(\ii)$ is a planar graph;
equivalently, no two inclined edges cross each other inside the strip.

\noindent (b) The boundary of any triangle or trapezoid formed by two 
consecutive inclined edges and horizontal segments between them is a 
directed cycle in $\Sigma_{i,j}(\ii)$. 
\end{theorem}

Our next goal is to compare the directed graphs $\Sigma(\ii)$ and 
$\Sigma({\ii'})$
when two reduced words $\ii$ and $\ii'$ are related by a $2$- or $3$-move. 
To do this, we associate to $\ii$ and $\ii'$ a permutation
$\sigma_{\ii',\ii}$ of $[1,m]$ defined as follows.
If $\ii$ and $\ii'$ are related by a trivial $2$-move in position $k$ then
$\sigma_{\ii',\ii} = (k-1,k)$, the transposition of $k-1$ and $k$;
if $\ii$ and $\ii'$ are related by a non-trivial $2$-move then
$\sigma_{\ii',\ii} = e$, the identity permutation of $[1,m]$;
finally, if $\ii$ and $\ii'$ are related by a $3$-move in position $k$ then
$\sigma_{\ii',\ii} = (k-2,k-1)$.
The following properties of $\sigma_{\ii',\ii}$ are immediate from the
definitions.

\begin{proposition}
\label{pr:sigma}
The permutation $\sigma_{\ii',\ii}$ sends $\ii$-bounded indices to 
$\ii'$-bounded ones. If the move that relates $\ii$ and $\ii'$ is 
non-trivial then its position $k$ is $\ii$-bounded, and 
$\sigma_{\ii',\ii}(k) = k$.
\end{proposition}

The relationship between the graphs $\Sigma(\ii)$ and $\Sigma(\ii')$ 
is now given as follows.

\begin{theorem}
\label{th:graph change}
Suppose two reduced words $\ii$ and $\ii'$ are related by a $2$- or $3$-move
in position $k$, and $\sigma = \sigma_{\ii',\ii}$ is the corresponding
permutation of $[1,m]$.
Let $a$ and $b$ be two distinct elements of $[1,m]$ such that at least 
one of them is $\ii$-bounded. Then 
\begin{equation}
\label{eq:preserving edges}
(a \to b) \in \Sigma(\ii) \Leftrightarrow (\sigma (a) \to \sigma (b)) \in 
\Sigma({\ii'}) \ ,
\end{equation}
with the following two exceptions.

\begin{enumerate}
\item  If the move that relates $\ii$ and $\ii'$ is non-trivial
then $(a \to k) \in \Sigma(\ii) \Leftrightarrow (k \to \sigma (a)) \in 
\Sigma({\ii'})$.

\item  If the move that relates $\ii$ and $\ii'$ is non-trivial, and 
$a \to k \to b$ in $\Sigma(\ii)$ then 
$\{a,b\} \in \Sigma(\ii) \Leftrightarrow \{\sigma (a), \sigma (b)\} 
\notin \Sigma({\ii'})$; 
furthermore, the edge $\{a,b\}\in\Sigma(\ii)$ can only be directed as $b\to a$.
\end{enumerate}
\end{theorem}

The following example illustrates the above results.  

\begin{example}
\label{ex:graph}
{\rm Let $\Pi$ be the Dynkin graph $A_4$, i.e., the chain formed by vertices $1, 2, 3$, and $4$.
Let $u = s_4 s_2 s_1 s_2 s_3 s_2 s_4 s_1$ and $v = s_2 s_1 s_3 s_2 s_4 s_1 s_3 s_2 s_1$
(in the standard realization of $W$ as the symmetric group $S_5$, 
with the generators $s_i = (i, i+1)$ (adjacent transpositions), 
the permutations $u$ and $v$ can be written in the one-line notation as 
$u=53241$ and $v=54312$). 
The graph $\Sigma (\ii)$ corresponding to the
reduced word $\ii=(2,1,-4,-2,-1,3,-2,2,-3,-2,4,1,-4,-1,3,2,1)$ of $(u,v)$
is shown on Fig.~\ref{fig:graph}. 
Here white (resp. black) vertices of each horizontal level $i$ correspond to entries of $\ii$ 
that are equal to $-i$ (resp. to $i$). 
Horizontal edges are shown by solid lines, inclined edges of 
type (ii) in Definition~\ref{def:edges} by dashed lines, and 
inclined edges of type (iii) by dotted lines.

\begin{figure}
\vskip 10pt
\centerline{\hbox{\epsfxsize=10cm\epsfbox{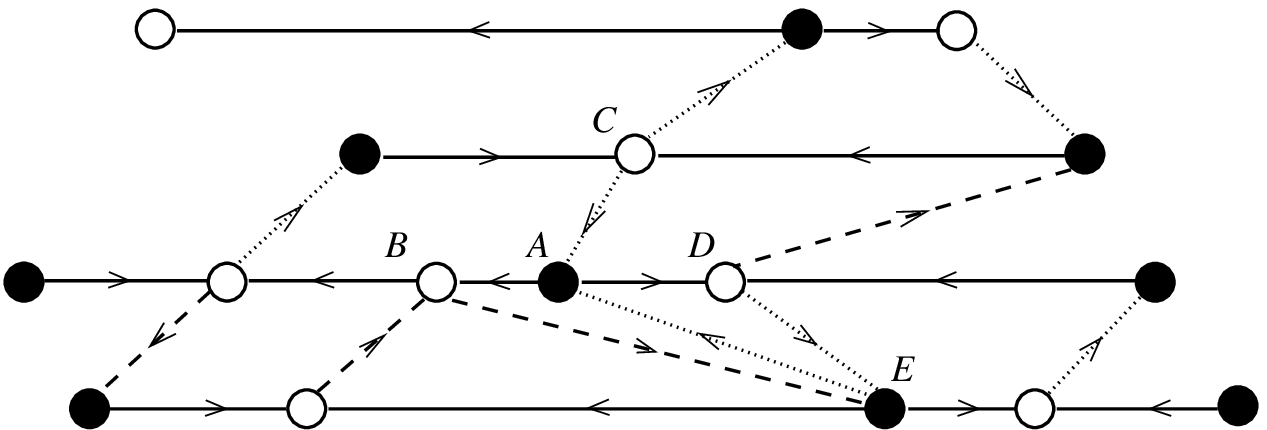}}}
\vskip 5pt
\caption[]{\label{fig:graph}
Graph $\Sigma(\ii)$ for 
type $A_4$. } 
\end{figure}

\begin{figure}[t]
\vskip 10pt
\centerline{\hbox{\epsfxsize=10cm\epsfbox{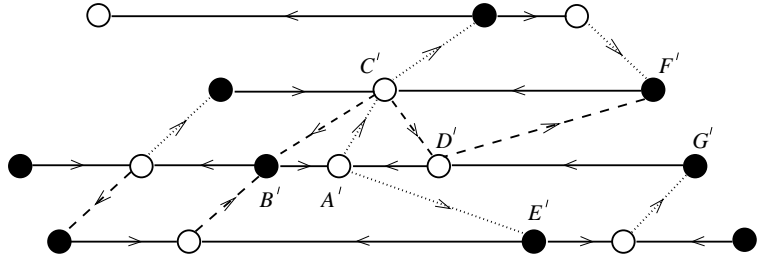}}}
\vskip 5pt
\caption[]{\label{fig:2-move}
Graph transformation under a non-trivial 2-move. } 
\end{figure}

\begin{figure}[b]
\vskip 10pt
\centerline{\hbox{\epsfxsize=10cm\epsfbox{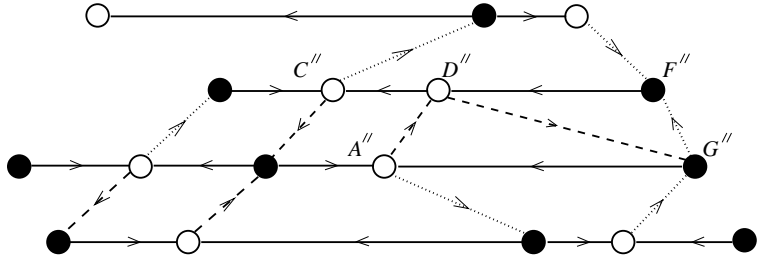}}}
\vskip 5pt
\caption[]{\label{fig:3-move}
Graph transformation under a 3-move. } 
\end{figure}

Now let $\ii'$ be obtained from $\ii$ by the (non-trivial) 2-move in position 8,
i.e., by interchanging $i_7 = -2$ with $i_8 = 2$. 
The corresponding graph $\Sigma (\ii')$ is shown on Fig.~\ref{fig:2-move}.

Notice that the edges of $\Sigma(\ii)$ that fall into 
the first exceptional case in Theorem~\ref{th:graph change} are 
$A\to B$, $C\to A$, and $A\to D$; by reversing their orientation,
one obtains the edges $B'\to A'$, $A'\to C'$, and $D'\to A'$ of $\Sigma(\ii')$.
The second exceptional case in Theorem~\ref{th:graph change} 
applies to two edges $B \to E$ and $D\to E$ of $\Sigma(\ii)$ 
and two ``non-edges" $\{C,B\}$ and $\{C,D\}$; the corresponding edges and 
non-edges of $\Sigma(\ii')$ are $C'\to B'$, $C'\to D'$, $\{E',B'\}$, and $\{E',D'\}$. 

Finally, consider the reduced word $\ii''$ obtained from $\ii'$ by the 3-move
in position 10, i.e., by replacing $(i'_8, i'_9, i'_{10}) = (-2,-3,-2)$ 
with $(-3, -2, -3)$. 
The corresponding graph $\Sigma (\ii'')$ is shown on Fig.~\ref{fig:3-move}.

Now the first exceptional case in Theorem~\ref{th:graph change} covers the edges 
$D'\to A'$, $C'\to D'$, $D'\to F'$, and $G'\to D'$ of $\Sigma(\ii')$, 
and the corresponding edges $A''\to D''$, $D''\to C''$, $F''\to D''$, and $D''\to G''$ 
of $\Sigma(\ii'')$.
The second exceptional case covers the edges $F'\to C'$ and $A'\to C'$, 
and non-edges $\{G',F'\}$ and $\{G',A'\}$ of $\Sigma(\ii'')$; 
the corresponding edges and non-edges of $\Sigma(\ii'')$ are
$G''\to F''$, $G''\to A''$, $\{C'',F''\}$, and $\{A'',C''\}$.}
\end{example}

\subsection{The groups $\Gamma_\ii$ and conjugacy theorems}
\label{sec:conjugacy}
As before, let $\ii = (i_1, \dots, i_m)$ be a reduced word for a pair $(u,v)$
of elements in a simply-laced Coxeter group $W$. 
By the general construction in Section~\ref{sec:symplectic transvections general}, the pair 
$(\Sigma(\ii), B(\ii))$
gives rise to a skew symmetric form $\Omega_{\Sigma(\ii)}$ on $\ZZ^m$,
and to a subgroup $\Gamma_{\Sigma (\ii), B(\ii)} \subset GL_m (\ZZ)$ 
generated by symplectic transvections.
We denote these symplectic transvections by $\tau_{k,\ii}$, and also 
abbreviate $\Omega_\ii = \Omega_{\Sigma(\ii)}$, and 
$\Gamma_\ii = \Gamma_{\Sigma (\ii), B(\ii)}$. 

\begin{theorem}
\label{th:conjugacy}
For any two reduced words $\ii$ and $\ii'$ for the same pair
$(u,v) \in W \times W$, the groups $\Gamma_\ii$ and $\Gamma_{\ii'}$
are conjugate to each other inside $GL_{m}(\ZZ)$.
\end{theorem}

Our proof of Theorem~\ref{th:conjugacy} is constructive.
In view of the Tits result quoted in Section~\ref{sec:Coxeter groups general},
it is enough to prove Theorem~\ref{th:conjugacy} in the case when 
$\ii$ and $\ii'$ are related by a 2- or 3-move. 
We shall construct the corresponding conjugating linear transformations 
explicitly.  To do this, let us define two linear maps
$\varphi^\pm_{\ii',\ii}: \ZZ^{m} \to \ZZ^{m}$.
For $v  \in \ZZ^{m}$, the vectors $\varphi^+_{\ii',\ii} (v) = v^+$ and
$\varphi^-_{\ii',\ii} (v) = v^-$ are defined as follows.
If $\ii$ and $\ii'$ are related by a trivial $2$-move and $l$ is
arbitrary, or if $\ii$ and $\ii'$ are related by a non-trivial move in 
position $k$ and $l \neq k$, then we set
\begin{equation}
\label{eq:phi non-critical}
\xi_l (v^+)  = \xi_l (v^-)  = \xi_{\sigma_{\ii',\ii} (l)} (v) \ ;
\end{equation}
for $l = k$ in the case of a non-trivial move, we set
\begin{equation}
\label{eq:phi critical}
\xi_k (v^+)  = \sum_{(a \to k) \in \Sigma(\ii)} \xi_a (v) - \xi_k (v) \ ; \quad
\xi_k (v^-)  = \sum_{(k \to b) \in \Sigma(\ii)} \xi_b  (v) - \xi_k  (v) \ .
\end{equation}

\begin{theorem}
\label{th:groupoid}
If two reduced words $\ii$ and $\ii'$ for the same pair
$(u,v) \in W \times W$ are related by a $2$- or $3$-move then
the corresponding linear maps $\varphi^+_{\ii',\ii}$ and 
$\varphi^-_{\ii',\ii}$ are invertible, and 
\begin{equation}
\label{eq:conjugacy}
\Gamma_{\ii'} = \varphi^+_{\ii',\ii} \circ \Gamma_{\ii} \circ 
(\varphi^+_{\ii',\ii})^{-1} 
= \varphi^-_{\ii',\ii} \circ \Gamma_{\ii} \circ (\varphi^-_{\ii',\ii})^{-1}\ .
\end{equation}
\end{theorem}

Our proof of Theorem~\ref{th:groupoid} is based on the following properties 
of the maps $\varphi^\pm_{\ii',\ii}$, which might be of independent interest. 

\begin{theorem}
\label{th:phi^2}

(a) The linear maps $\varphi^\pm_{\ii',\ii}$ satisfy:
\begin{equation}
\label{eq:phi inverses}
\varphi^-_{\ii,\ii'} \circ \varphi^+_{\ii',\ii} =
\varphi^+_{\ii,\ii'} \circ \varphi^-_{\ii',\ii} = {\rm Id} \ .
\end{equation}

\noindent (b) If the move that relates $\ii$ and $\ii'$ is non-trivial 
in position $k$ then
\begin{equation}
\label{eq: tau = phi^2}
\varphi^+_{\ii,\ii'} \circ \varphi^+_{\ii',\ii} = \tau_{k,\ii} \ .
\end{equation}

\noindent (c) For any $\ii$-bounded index $l \in [1,m]$, we have 
\begin{equation}
\label{eq:phi-intertwiner}
\varphi^+_{\ii',\ii} \circ \tau_{l,\ii} = 
\tau_{\sigma_{\ii',\ii} (l),\ii'} \circ \varphi^+_{\ii',\ii} 
\end{equation}
unless the move that relates $\ii$ and $\ii'$ is non-trivial in position $k$, 
and $(l \to k) \in \Sigma_\ii$. 

\end{theorem}

\subsection{Enumerating $\Gamma_{\Sigma, B}( \FF_2)$-orbits in $\FF_2^\Sigma$}
\label{sec:orbit enumeration}
Let $\Sigma$ and $B$ have the same meaning as in 
Section~\ref{sec:symplectic transvections general}, 
and let $\Gamma = \Gamma_{\Sigma, B}( \FF_2)$ be the corresponding group of 
linear transformations of the vector space $\FF_2^\Sigma$. 

The following definition is motivated by the results in 
\cite{Janssen,ssv1,ssv2}. 

\begin{figure}
\vskip 10pt
\centerline{\hbox{\epsfxsize=6cm\epsfbox{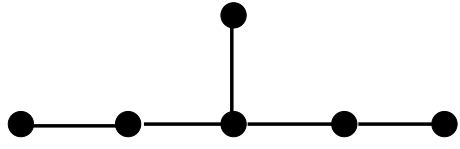}}}
\vskip 5pt
\caption[]{\label{fig:E6}
Dynkin graph $E_6$. } 
\end{figure}

\begin{definition}
\label{def:E6-compatible graphs}
{\rm A finite (non-directed) graph is $E_6$-compatible if it is connected, 
and it contains an induced subgraph with $6$ vertices isomorphic to the Dynkin 
graph $E_6$ (see Fig.~\ref{fig:E6}).}
\end{definition}

\begin{theorem}
\label{th:number of orbits}
Suppose that the induced subgraph of $\Sigma$ with the set of vertices 
$B$ is $E_6$-compatible. 
Then the number of $\Gamma$-orbits in $\FF_2^\Sigma$ is equal to 
$$2^{\#(\Sigma \setminus B)} \cdot (2 + 2^{\dim (\FF_2^B  \cap 
{\rm Ker} \ \overline \Omega)}) \ ,$$
where $\overline \Omega$ denotes the $\FF_2$-valued bilinear form on 
$\FF_2^\Sigma$ obtained by reduction modulo $2$ from the form 
$\Omega = \Omega_\Sigma$ in (\ref{eq:Omega}). 
\end{theorem}

Theorem~\ref{th:number of orbits} has the following corollary which 
generalizes the main enumeration result in \cite{ssv1,ssv2}. 

\begin{corollary}
\label{cor:number of components}
Let $u$ and $v$ be two elements of a simply-laced Coxeter group $W$, and 
suppose that for some reduced word $\ii \in R(u,v)$, the induced subgraph 
of $\Sigma(\ii)$ with the set of vertices $B(\ii)$ is $E_6$-compatible. 
Then the number of $\Gamma_\ii (\FF_2)$-orbits in $\FF_2^{m}$ is equal to 
$3 \cdot 2^s$, where $s$ is the number of indices $i \in \Pi$ such that 
some (equivalently, any) reduced word for $(u,v)$ has an entry $\pm i$.
\end{corollary}

\section{Proofs of results in Section~\ref{sec:graphs}}
\label{sec:proofs-graphs}

\subsection{Proof of Proposition~\ref{pr:boundary vertex}}
\label{sec:ProofBoundary}
By the definition of $\ii$-bounded indices, we have $k^- \in [1,m]$ for any 
$k \in S$. 
Now pick $b \in S$ with the smallest value of $b^-$, and set $a = b^-$. 
Clearly, $a \notin S$, and $\{a,b\}$ is a horizontal edge in $\Sigma (\ii)$. 
We claim that $b$ is the only vertex in $S$ such that 
$\{a,b\} \in \Sigma (\ii)$. 
Indeed, if $\{a,c\} \in \Sigma (\ii)$ for some $c \neq b$
then $c^- < a$, in view of Definition~\ref{def:edges}.
Because of the way $b$ was chosen, we have $c \notin S$, as required.  
\endproof

\subsection{Proof of Theorem~\ref{th:strip}}
\label{sec:ProofStrips}
In the course of the proof, we fix a reduced word $\ii \in R(u,v)$, 
and an edge $\{i,j\} \in \Pi$;
we shall refer to the $(i,j)$-strip of $\Sigma(\ii)$ as simply the strip. 
For any vertex $A = A_k = (k,y)$ in the strip, we set $y(A) = y$, and 
$\veps (A) = \veps (i_k)$;
we call $y(A)$ the \emph{level}, and $\veps (A)$ the \emph{sign} of $A$. 
We also set
$$c(A) = y(A) \veps(A) \ ,$$
and call $c(A)$ the \emph{charge} of a vertex $A$.
Finally, we linearly order the vertices by setting $A_k \prec A_l$ if $k < l$,
 i.e., if the vertex $A_k$ is to the left of $A_l$. 
In these terms, one can describe inclined edges in the strip as follows.

\begin{lemma}
\label{lem:inclined edges}
A vertex $B$ is the left end of an inclined edge in the strip  
if and only if it satisfies the following two conditions:

\noindent (1) $B$ is not the leftmost vertex in the strip, and the 
preceding vertex $A$  has opposite charge $c(A) = - c(B)$. 

\noindent (2) there is a vertex $C$ of opposite level $y(C) = - y(B)$ 
that lies to the right of $B$.

Under these conditions, an inclined edge with the left end $B$ is unique,   
and its right end is the leftmost vertex $C$ satisfying (2). 
\end{lemma}

This is just a reformulation of conditions (ii) and (iii) in 
Definition~\ref{def:edges}. 
\endproof

\begin{lemma}
\label{lem:inclined edges 2} 
Suppose $A \prec C \prec C'$ are three vertices in the strip such that 
$c(A) = - c(C)$, and $y(C) = - y(C')$. 
Then there exists a vertex $B$ such that $A \prec B \preceq C$, and $B$ 
is the left end of an inclined edge in the strip. 
\end{lemma}

\proof 
Let $B$ be the leftmost vertex such that $A \prec B \preceq C$ and 
$c(B) = - c(A)$. 
Clearly, $B$ satisfies condition (1) in Lemma~\ref{lem:inclined edges}. 
It remains to show that $B$ also satisfies condition (2); that is, 
we need to find a vertex of opposite level to $B$ that lies to the 
right of $B$. Depending on the level of $B$, either $C$ or $C'$ is such 
a vertex, and we are done.
\endproof

Now everything is ready for the proof of Theorem~\ref{th:strip}. 
To prove part (a), assume that $\{B,C\}$ and $\{B',C'\}$ are two 
inclined edges that cross each other inside the strip.
Without loss of generality, assume that $B \prec C$, $B' \prec C'$, 
and $C \prec C'$. Then we must have $B' \prec C$ 
(otherwise, our inclined edges would not cross). 
Since $\{B',C'\}$ is an inclined edge, and $B' \prec C \prec C'$, 
Lemma~\ref{lem:inclined edges} implies that $y(C) = y(B')$. 
Therefore, $y(B) = - y(C) = - y(B')$. 
Again applying Lemma~\ref{lem:inclined edges} to the inclined edge $\{B',C'\}$,
we conclude that $B \prec B'$, i.e., we must have 
$B \prec B' \prec C \prec C'$.  
But then, by the same lemma, $\{B,C\}$ cannot be an inclined edge,
providing a desired contradiction. 

To prove part (b), consider two consecutive inclined edges $\{B,C\}$ 
and $\{B',C'\}$. Again we can assume without loss of generality that 
$B \prec C$, $B' \prec C'$, and $C \prec C'$.
Let $P$ be the boundary of the polygon with vertices $B, C, B'$, and $C'$. 
By Lemma~\ref{lem:inclined edges}, the leftmost vertex of $P$ is $B$, 
the rightmost vertex is $C'$, and $P$ does not contain a vertex $D$ such that
$B' \prec D \prec C'$; in particular, we have either $C \preceq B'$ or 
$C = C'$. Now we make the following crucial observation: all 
the vertices $D$ on $P$ such that 
$B \prec D \prec B'$ must have the same charge $c(D) = c(B)$.
Indeed, assume that $c(D) = - c(B)$ for some $D$ with $B \prec D \prec B'$.
Then Lemma~\ref{lem:inclined edges 2} implies that some $B''$ with 
$B \prec B'' \preceq D$ 
is the left end of an inclined edge; but this contradicts our assumption that 
$\{B,C\}$ and $\{B',C'\}$ are two \emph{consecutive} inclined edges. 
We see that $c(D) = c(B)$ for any vertex $D \in P \setminus \{B',C'\}$.
Combining this fact with condition (1) in Lemma~\ref{lem:inclined edges} 
applied to the inclined edge $\{B',C'\}$ with the left end $B'$, we conclude 
that $c(B') = - c(B)$. 
Remembering the definition of charge, the above statements can be 
reformulated as follows: 
$B'$ has the same (resp. opposite) sign with all vertices of 
opposite (resp. same) level in $P \setminus \{C'\}$.   
Using the definition of directions of edges in Definition~\ref{def:edges},
we obtain:
 
\noindent 1. Horizontal edges on opposite sides of $P$ are directed 
opposite way since their left ends have opposite signs. 
 
\noindent 2. Suppose $B'$ is the right end of 
a horizontal edge $\{A,B'\}$ in $P$.
Then exactly one of the edges $\{A,B'\}$  and
$\{B',C'\}$ is directed towards $B'$ since their left ends $A$ and $B'$ 
have opposite signs. 
 
\noindent 3. The same argument shows that if $C'$ is the right end of 
a horizontal edge $\{A,C'\}$ in $P$ 
then exactly one of the edges $\{A,C'\}$  and
$\{B',C'\}$ is directed towards $C'$. 
 
\noindent 4. Finally, if $B$ is the left end of 
a horizontal edge $\{B,D\}$ in $P$ then 
exactly one of the edges $\{B,C\}$ and
$\{B,D\}$ is directed towards $B$. 
 
These facts imply that $P$ is a directed cycle, which completes the proof of 
Theorem~\ref{th:strip}.
\endproof

\subsection{Proof of Theorem~\ref{th:graph change}}
\label{sec:ProofGraphChange}
Let us call a pair of indices $\{a,b\}$ \emph{exceptional} 
(for $\ii$ and $\ii'$) 
if it violates (\ref{eq:preserving edges}).
We need to show that exceptional pairs are precisely those in two 
exceptional cases in Theorem~\ref{th:graph change}; to do this, 
we shall examine the relationship between the 
corresponding strips in $\Sigma(\ii)$ and $\Sigma(\ii')$. 
Let us consider the following three cases:

\noindent {\bf Case 1 (trivial $2$-move).} Suppose $i_k = i'_{k-1} = i_0$, 
$i_{k-1} = i'_{k} = j_0$, and $i_l = i'_l$ for $l \notin \{k-1,k\}$, where 
$i_0, j_0 \in \tilde \Pi$ are such that 
$|i_0| \neq |j_0|$ and $\{i_0,j_0\} \notin \tilde \Pi$. 

If both $i$ and $j$ are different from $|i_0|$ and $|j_0|$ then the 
strip $\Sigma_{i,j}(\ii)$ is identical to $\Sigma_{i,j}(\ii')$, 
and so does not contain exceptional pairs.  

If say $i = |i_0|$ but $j \neq |j_0|$ then the only vertex in 
$\Sigma_{i,j}(\ii)$ but not in
$\Sigma_{i,j}(\ii')$ is $A_k$ (in the notation of 
Section~\ref{sec:ProofStrips}), while the only vertex in 
$\Sigma_{i,j}(\ii')$ but not in
$\Sigma_{i,j}(\ii)$ is $A'_{k-1} = A'_{\sigma(k)}$. 
The vertex $A_k$ has the same level and sign and so the same charge as 
the vertex $A'_{\sigma(k)}$ in $\Sigma_{i,j}(\ii')$; 
by Lemma~\ref{lem:inclined edges}, there are no exceptional pairs in the strip 
$\Sigma_{i,j}(\ii)$.

Finally, suppose that $\{i,j\} = \{|i_0|, |j_0|\}$; in particular, 
in this case we have
$\{|i_0|, |j_0|\} \in \Pi$, hence $\veps (i_0) = - \veps (j_0)$. 
Now the only vertices in $\Sigma_{i,j}(\ii)$ but not in
$\Sigma_{i,j}(\ii')$ are $A_k$ and $A_{k-1}$, while the only vertices in 
$\Sigma_{i,j}(\ii')$ but not in $\Sigma_{i,j}(\ii)$ are 
$A'_{k-1} = A'_{\sigma(k)}$ and $A'_{k} = A'_{\sigma(k-1)}$. 
Since $A_k$ and $A_{k-1}$ are of opposite level and opposite sign, 
they have the same charge,
which is also equal to the charge of $A'_{\sigma(k-1)}$ and $A'_{\sigma(k)}$. 
Again using Lemma~\ref{lem:inclined edges}, we see that the strip in 
question also does not contain exceptional pairs. 

\noindent {\bf Case 2 (non-trivial $2$-move).} Suppose 
$i_k = i'_{k-1} = i_0 \in \tilde \Pi$, 
$i_{k-1} = i'_{k} = - i_0$, and $i_l = i'_l$ for $l \notin \{k-1,k\}$. 
Interchanging if necessary $\ii$ and $\ii'$, we can and will assume that 
$i_0 \in \Pi$. Clearly, an exceptional pair can only belong to an 
$(i,j)$-strip with $i = i_0$. In our case, the location of all vertices in 
$\Sigma_{i,j}(\ii)$ and $\Sigma_{i,j}(\ii')$
is the same; the only difference between the two strips is that the 
vertices $A_{k-1}$ and $A_k$ in $\Sigma_{i,j}(\ii)$
have opposite signs and hence opposite charges to their counterparts in 
$\Sigma_{i,j}(\ii')$. 
It follows that exceptional pairs of vertices of the same level are precisely 
horizontal edges containing $A_k$, i.e., $\{A_{k-1}, A_k\}$ and $\{A_k, C\}$,
where $C$ is the right neighbor of $A_k$ of the same level (note that $C$ 
does not necessarily exist). 
Since $\veps (i_k) = \veps (i'_{k-1}) = +1$, and 
$\veps (i_{k-1}) = \veps (i'_{k}) = -1$, we have 
$$(A_k \to A_{k-1}) \in \Sigma(\ii), \,\, (A_k \to C) \in \Sigma(\ii),$$
$$(A'_{k-1} \to A'_{k}) \in \Sigma(\ii'), \,\, (C' \to A'_k) \in\Sigma(\ii'),$$
so both pairs $\{A_{k-1}, A_k\}$ and $\{A_k, C\}$ fall into the first 
exceptional case in Theorem~\ref{th:graph change}.

Let us now describe exceptional pairs corresponding to inclined edges. 
Let $B$ be the vertex of the opposite level to $A_k$ and closest to $A_k$ 
from the right (as the vertex $C$ above, $B$ does not necessarily exist). 
By Lemma~\ref{lem:inclined edges}, the left end of an exceptional 
inclined pair can only be $A_{k-1}$, $A_k$, or the leftmost of $B$ and $C$; 
furthermore, the corresponding inclined edges can only be $\{A_{k-1},B\}$,  
$\{A_k,B\}$, or $\{B,C\}$. 
We claim that all these three pairs are indeed exceptional, and each of 
them falls into one of the exceptional cases in Theorem~\ref{th:graph change}. 

Let us start with $\{B,C\}$.
Since $A_k$ is the preceding vertex to the leftmost member of $\{B,C\}$, and 
it has opposite charges in the two strips, Lemma~\ref{lem:inclined edges}
implies that $\{B,C\}$ is an edge in precisely one of the strips.
By Theorem~\ref{th:strip} (b), the triangle with vertices $A_k$, $B$, and $C$
is a directed cycle in the corresponding strip. 
Thus the pair $\{B,C\}$ falls into the second exceptional case in 
Theorem~\ref{th:graph change}.

The same argument shows that $\{A_{k-1},B\}$ falls into the second 
exceptional case in Theorem~\ref{th:graph change} provided one of 
$A_{k-1}$ and $B$ is $\ii$-bounded, i.e., 
$A_{k-1}$ is not the leftmost vertex in the strip. 
As for $\{A_k,B\}$, it is an edge in both strips, and it has opposite 
directions in them because its left end $A_k$ has opposite signs there.  
Thus $\{A_k,B\}$ falls into the first exceptional case in 
Theorem~\ref{th:graph change}. 

It remains to show that the exceptional pairs (horizontal and inclined) just 
discussed exhaust all possibilities for the 
two exceptional cases in Theorem~\ref{th:graph change}.
This is clear because by the above analysis, the only possible edges through 
$A_k$ in $\Sigma (\ii)$ are $(A_k \to A_{k-1})$, $(A_k \to C)$, and 
$(B \to A_k)$ with $B$ of the kind described above. 

\noindent {\bf Case 3 ($3$-move).} Suppose $i_k = i_{k-2} = i'_{k-1} = i_0$, 
$i_{k-1} = i'_{k} = i'_{k-2} = j_0$ for some $\{i_0, j_0\} \in \Pi$, and 
$i_l = i'_l$ for $l \notin \{k-2, k-1,k\}$ (the case when 
$\{i_0, j_0\} \in -\Pi$ is totally similar).
As in the previous case, we need to describe all exceptional pairs.  

First an exceptional pair can only belong to an $(i,j)$-strip with 
at least one of $i$ and $j$ equal to $i_0$ or $j_0$. 
Next let us compare the $(i_0,j_0)$-strips in $\Sigma(\ii)$ and 
$\Sigma(\ii')$. The location of all vertices in these two strips 
is the same with the exception of $A_{k-2}, A_{k-1}$, and $A_k$ in the 
former strip, and their counterparts $A'_{k-2} = A'_{\sigma(k-1)}, 
A'_{k-1} = A'_{\sigma(k-2)}$, and $A'_k$ in the latter strip. 
Each of the six exceptional vertices has sign $+ 1$; so its level is equal 
to its charge. These charges (or levels) are given as follows:
$$c(A_{k-2}) = c(A'_{\sigma (k-2)}) = c(A_{k}) = -1, \,\, 
c(A_{k-1}) = c(A'_{\sigma(k-1)}) = c(A'_{k}) = 1 \ .$$
Let $B$ (resp. $B'$) denote the vertex in both strips which is the closest 
from the right to $A_k$ on the same (resp. opposite) level; note that $B$ 
or $B'$ may not exist. Since the trapezoid $T$ with vertices 
$A_{k-2}, A_{k-1}, B'$, and $B$ in $\Sigma_{i_0,j_0}(\ii)$ 
is in the same relative position to all outside vertices as 
the trapezoid $T'$ with vertices $A'_{\sigma(k-2)}, A'_{\sigma(k-1)}, B'$, a
nd $B$ in $\Sigma_{i_0,j_0}(\ii')$, it follows that every exceptional pair 
is contained in $T$. 
An inspection using Lemma~\ref{lem:inclined edges} shows that $T$ contains 
the directed edges
$$A_{k-2} \to A_k \to A_{k-1} \to B' \to A_k \to B$$
and does not contain any of the edges $\{A_{k-2}, B\}$, $\{A_{k-2}, B'\}$, 
or $\{A_{k-1}, B\}$. Similarly (or by interchanging $\ii$ and $\ii'$), 
we conclude that $T'$ contains the directed edges
$$A'_{\sigma(k-1)} \to A'_k \to A'_{\sigma(k-2)} \to B \to A'_k \to B'$$
and does not contain any of the edges $\{A'_{\sigma(k-1)}, B'\}$, 
$\{A'_{\sigma(k-1)}, B'\}$, or $\{A'_{\sigma(k-2)}, B'\}$. 
Furthermore, $\{B,B'\}$ is an edge in precisely one of the strips 
(since the preceding vertices
$A_k$ and $A'_k$ have opposite charges); and precisely one of the pairs 
$\{A_{k-2}, A_{k-1}\}$
and $\{A'_{\sigma (k-1)}, A'_{\sigma (k-2)}\}$ 
is an edge in its strip provided 
$A_{k-2}$ is not the leftmost vertex (since their left ends $A_{k-2}$ and 
$A'_{\sigma (k-1)}$ have opposite charges). 

Comparing this information for the two trapezoids, we see that the exceptional
pairs in $T$ are all pairs of vertices in $T$ with the exception of two 
diagonals
$\{A_{k-2}, B'\}$ and $\{A_{k-1}, B\}$ (and also of $\{A_{k-2}, A_{k-1}\}$ if 
$A_{k-2}$ is the leftmost vertex in the strip).
By inspection based on Theorem~\ref{th:strip} (b), all these exceptional 
pairs fall 
into the two exceptional cases in Theorem~\ref{th:graph change}.

A similar (but much simpler) analysis shows that any $(i,j)$-strip
with precisely one of $i$ and $j$ belonging to $\{i_0,j_0\}$ does not 
contain extra exceptional pairs, and also has no inclined edges through 
$A_k$ or $A'_k$. We conclude that all the exceptional pairs are contained 
in the above trapezoid $T$. 
The fact that these exceptional pairs exhaust all possibilities for the 
two exceptional cases in Theorem~\ref{th:graph change} 
is clear because by the above analysis, the only edges through $A_k$ in 
$\Sigma (\ii)$ are those connecting $A_k$ with the vertices of $T$. 
Theorem~\ref{th:graph change} is proved. 
\endproof

\section{Proofs of results in Section~\ref{sec:conjugacy}}
\label{sec:proofs-abstract}

We have already noticed that Theorem~\ref{th:conjugacy} follows from 
Theorem~\ref{th:groupoid}. 
Let us first prove Theorem~\ref{th:phi^2} and then deduce 
Theorem~\ref{th:groupoid} from it.

\subsection{Proof of Theorem~\ref{th:phi^2}}
We fix reduced words $\ii$ and $\ii'$ related by a 2- or 3-move, and 
abbreviate $\sigma = \sigma_{\ii',\ii} = \sigma_{\ii,\ii'}$ and 
$\varphi^+ = \varphi^+_{\ii',\ii}$. Let us first prove parts (a) and (b).  
We shall only prove the first equality in (\ref{eq:phi inverses}); the 
proof of the second one and of (\ref{eq: tau = phi^2}) is completely similar. 
Let $v \in \ZZ^m$, $v^+ = \varphi^+ (v)$, and 
$v' = \varphi^-_{\ii,\ii'} (v^+)$;
thus we need to show that $v = v'$, i.e., that $\xi_l (v) = \xi_l (v')$ 
for all $l \in [1,m]$.  Note that the permutation $\sigma$ is an involution.
In view of (\ref{eq:phi non-critical}), this implies the desired equality 
$\xi_l (v) = \xi_l (v')$ in all the cases except the following one: 
the move that relates $\ii$ and $\ii'$ is non-trivial in position $k$, and 
$l = k$. To deal with this case, we use the first exceptional case in 
Theorem~\ref{th:graph change} which we can write as
$$(k \to b) \in \Sigma(\ii') \Leftrightarrow (\sigma (b) \to k) \in 
\Sigma({\ii}) \ .$$
Combining this with the definitions (\ref{eq:phi non-critical}) and 
(\ref{eq:phi critical}), we obtain
$$\xi_k (v') = \sum_{(k \to b) \in \Sigma(\ii')} \xi_b  (v^+) - \xi_k  (v^+)$$
$$ = \sum_{(\sigma (b) \to k) \in \Sigma(\ii)} \xi_{\sigma (b)}(v) -
(\sum_{(a \to k) \in \Sigma(\ii)} \xi_{a}(v) - \xi_k (v)) = \xi_k (v) \ ,$$
as required. 

We deduce part (c) from the following lemma which says that
the maps $(\varphi^\pm_{\ii',\ii})^*$ induced by
$\varphi^\pm_{\ii',\ii}$ ``almost" transform the form $\Omega_{\ii'}$ into
$\Omega_\ii$. 

\begin{lemma}
\label{lem:Omega-transform}
If the move that relates $\ii$ and $\ii'$ is trivial then
$$(\varphi^+_{\ii',\ii})^* (\Omega_{\ii'}) =
(\varphi^-_{\ii',\ii})^* (\Omega_{\ii'}) = \Omega_{\ii} \ .$$
If the move that relates $\ii$ and $\ii'$ is non-trivial in position $k$ then
\begin{equation}
\label{eq:Omega-transform}
(\varphi^+_{\ii',\ii})^* (\Omega_{\ii'}) =
(\varphi^-_{\ii',\ii})^* (\Omega_{\ii'}) = \Omega_{\ii} - 
\doublesubscript{\sum}{(a \to k \to b) \in \Sigma(\ii)}{a,b \notin B(\ii)}
\xi_a \wedge \xi_b \ .
\end{equation}
\end{lemma}

\proof
We will only deal with  
$(\varphi^+)^* (\Omega_{\ii'})=(\varphi^+_{\ii',\ii})^* (\Omega_{\ii'})$; 
the form $(\varphi^-_{\ii',\ii})^*(\Omega_{\ii'})$ can be treated in the 
same way. 
By the definition,
$$(\varphi^+)^* (\Omega_{\ii'}) =  
\sum_{(a' \to b') \in \Sigma(\ii')} (\varphi^+)^* \xi_{a'} \wedge 
(\varphi^+)^* \xi_{b'} \ .$$
The forms $(\varphi^+)^* \xi_{a'}$ are given by (\ref{eq:phi non-critical}) 
and (\ref{eq:phi critical}). 
In particular, if $\ii$ and $\ii'$ are related by a trivial move then 
$(\varphi^+)^* \xi_{a'} = \xi_{\sigma(a')}$ for any $a' \in [1,m]$; 
by Theorem~\ref{th:graph change}, in this case we have
$$(\varphi^+)^* (\Omega_{\ii'}) = \sum_{(a \to b) \in \Sigma(\ii)} \xi_{a} 
\wedge \xi_{b} \,$$ 
as claimed.

Now suppose that $\ii$ and $\ii'$ are related by a non-trivial move in 
position $k$. Then we have 
\begin{eqnarray*}
(\varphi^+)^* (\Omega_{\ii'}) &=& 
\doublesubscript{\sum}{(\sigma(a) \to \sigma(b)) \in \Sigma(\ii')}{a,b \neq k} 
\xi_a \wedge \xi_b \\
&+& \sum_{(k \to \sigma (a')) \in \Sigma(\ii')} (\sum_{(a \to k) \in 
\Sigma(\ii)} \xi_{a} - \xi_k) \wedge \xi_{a'} \\
&+& \sum_{(\sigma (b) \to k) \in \Sigma(\ii')} \xi_{b} \wedge 
(\sum_{(a \to k) \in \Sigma(\ii)} \xi_{a} - \xi_k) \ .
\end{eqnarray*}
Using the second exceptional case in Theorem~\ref{th:graph change}, we can 
rewrite the first summand as 
$$\doublesubscript{\sum}{(a \to b) \in \Sigma(\ii)}{a,b \neq k} \xi_a \wedge \xi_b
+ \doublesubscript{\sum}{(a \to k \to b) \in \Sigma(\ii)}{\{a,b\} \cap B(\ii) \neq \emptyset} 
\xi_a \wedge \xi_b \ .$$
Similarly, using the first exceptional case in Theorem~\ref{th:graph change}, 
we can rewrite the last two summands as 
$$\sum_{(a \to k) \in \Sigma(\ii)} \xi_{a} \wedge \xi_k +  \sum_{(k \to b) 
\in \Sigma(\ii)} \xi_{k} \wedge \xi_b
- \sum_{(a \to k \to b) \in \Sigma(\ii)} \xi_a \wedge \xi_b$$
(note that the missing term 
$$\doublesubscript{\sum}{(a \to k) \in \Sigma(\ii)}{(a' \to k) \in 
\Sigma(\ii)} \xi_a \wedge \xi_{a'}$$
is equal to $0$). 
Adding up the last two sums, we obtain (\ref{eq:Omega-transform}).
\endproof

Now everything is ready for the proof of Theorem~\ref{th:phi^2} (c). 
Since $l$ is assumed to be $\ii$-bounded, Lemma~\ref{lem:Omega-transform} 
implies that 
$\Omega_{\ii} (v, e_l) = \Omega_{\ii'} (\varphi^+(v), \varphi^+(e_l))$
for any $v \in \ZZ^m$. On the other hand, since the case when the move 
that relates $\ii$ and $\ii'$ is non-trivial 
in position $k$, and $(l \to k) \in \Sigma_\ii$, is excluded, we have 
$\varphi^+(e_l) = \pm e_{\sigma(l)}$ (with the minus sign for $l = k$ only).
Therefore, our assumptions on $l$ imply that 
$$\Omega_{\ii} (v, e_l) \varphi^+ (e_l) =  
\Omega_{\ii'} (\varphi^+(v), e_{\sigma(l)}) e_{\sigma(l)} \ .$$
Remembering the definition (\ref{eq:transvections}) of symplectic 
transvections, we conclude that 
$$(\tau_{\sigma (l),\ii'} \circ \varphi^+)(v) = 
\varphi^+(v) - \Omega_{\ii'} (\varphi^+(v), e_{\sigma(l)}) e_{\sigma(l)}$$
$$= \varphi^+(v) - \Omega_{\ii} (v, e_l) \varphi^+ (e_l)
= (\varphi^+ \circ \tau_{l,\ii})(v) \ ,$$
as required.
This completes the proof of Theorem~\ref{th:phi^2}.
\endproof

\begin{remark} {\rm It is possible to modify all skew symmetric forms 
$\Omega_\ii$ 
without changing the corresponding groups $\Gamma_\ii$ in such a way that 
the modified forms will be preserved by the maps $(\varphi^\pm_{\ii',\ii})^*$. 
There are several ways to do it. 
Here is one ``canonical" solution: replace each $\Omega_\ii$ by the form
$$\tilde \Omega_\ii = \Omega_\ii - \frac{1}{2} \sum \veps(i_k) \xi_k 
\wedge \xi_l \ ,$$
where the sum is over all pairs of $\ii$-unbounded indices $k < l$ such that 
$\{|i_k|, |i_l|\} \in \Pi$.
It follows easily from Lemma~\ref{lem:Omega-transform} that 
$(\varphi^\pm_{\ii',\ii})^* (\tilde \Omega_{\ii'}) = \tilde \Omega_{\ii}$.  
Unfortunately, the forms $\tilde \Omega_\ii$ are not defined over $\ZZ$; 
in particular,  they cannot be reduced to bilinear forms over $\FF_2$.} 
\end{remark}

\subsection{Proof of Theorem~\ref{th:groupoid}}
The fact that $\varphi^+_{\ii',\ii}$ and $\varphi^-_{\ii',\ii}$ are invertible
follows from (\ref{eq:phi inverses}). 
To prove (\ref{eq:conjugacy}), it remains to show that 
$\varphi^+_{\ii',\ii} \circ \tau_{l,\ii} \circ (\varphi^+_{\ii',\ii})^{-1} 
\in \Gamma_{\ii'}$ for any $\ii$-bounded index $l \in [1,m]$. 
This follows from (\ref{eq:phi-intertwiner}) unless the move that relates 
$\ii$ and $\ii'$ is non-trivial in position $k$, and 
$(l \to k) \in \Sigma_\ii$. 
In this exceptional case, we conclude by interchanging $\ii$ and $\ii'$ in 
(\ref{eq:phi-intertwiner}) that 
$$\varphi^+_{\ii,\ii'} \circ \tau_{\sigma_{\ii',\ii} (l),\ii'}  = 
\tau_{l,\ii} \circ \varphi^+_{\ii,\ii'} \ .$$
Using (\ref{eq: tau = phi^2}), we obtain that  
$$\varphi^+_{\ii',\ii} \circ \tau_{l,\ii} \circ (\varphi^+_{\ii',\ii})^{-1} = 
(\varphi^+_{\ii',\ii} \circ \varphi^+_{\ii,\ii'}) \circ 
\tau_{\sigma_{\ii',\ii} (l),\ii'} 
\circ (\varphi^+_{\ii',\ii} \circ \varphi^+_{\ii,\ii'})^{-1} = \tau_{k,\ii'} 
\circ \tau_{\sigma_{\ii',\ii} (l),\ii'} 
\circ \tau_{k,\ii'}^{-1} \in \Gamma_{\ii'},$$
as required. 
This completes the proofs of Theorems~\ref{th:groupoid} and \ref{th:conjugacy}.
\endproof

\section{Proofs of results in Section~\ref{sec:orbit enumeration}}
\label{sec:proofs-enumeration}

\subsection{Description of $\Gamma$-orbits}
In this section we shall only work over the field $\FF_2$.
Therefore we find it convenient to change our notation a little bit. 
Let $V$ be a finite-dimensional vector space over $\FF_2$
with a skew-symmetric $\FF_2$-valued form $\Omega$
(i.e., $\Omega (v,v) = 0$ for any $v \in V$).
For any $v \in V$,  let $\tau_v: V \to V$ denote the corresponding symplectic
transvection acting by $\tau_v (x) = x - \Omega (x,v) v$.
Fix a linearly independent subset $B \subset V$, and let $\Gamma$ 
be the subgroup of $GL(V)$ generated by the transvections $\tau_b$ for 
$b \in B$.
We make $B$ the set of vertices of a graph with $\{b,b'\}$ an edge whenever
$\Omega (b,b') = 1$.

We shall deduce Theorem~\ref{th:number of orbits} from the following 
description of the $\Gamma$-orbits in $V$ in the case when the graph $B$ is 
$E_6$-compatible (see Definition~\ref{def:E6-compatible graphs}). 

Let $U \subset V$ be the linear span of $B$.
The group $\Gamma$ preserves each parallel translate $(v + U) \in V/U$ of
$U$ in $V$, so we only need to describe $\Gamma$-orbits in each $v + U$.

Let us first describe one-element orbits, i.e., $\Gamma$-fixed points in
each ``slice" $v + U$.
Let $V^\Gamma \subset V$ denote the subspace of $\Gamma$-invariant vectors, 
and $K \subset U$ denote the kernel of the restriction $\Omega|_U$.

\begin{proposition}
\label{pr:Gamma-invariants}
If $\Omega(K,v + U) = 0$ then $(v + U) \cap V^\Gamma$ is a
parallel translate of $K$; otherwise, this intersection is empty.
\end{proposition}

\proof
Suppose the intersection $(v + U) \cap V^\Gamma$ is non-empty; without loss of 
generality, we can assume that $v$ is $\Gamma$-invariant. 
By the definition, $v \in V^\Gamma$ if and only if $\Omega (u,v) = 0$ for all
$u \in U$.
In particular, $\Omega(K,v) = 0$, hence $\Omega(K,v + U) = 0$. 
Furthermore, an element $v + u$ of $v + U$ is $\Gamma$-invariant if and only if
$u \in K$, and we are done.
\endproof

Following \cite{Janssen}, we choose a function $Q: V \to \FF_2$ satisfying 
the following properties:
\begin{equation}
\label{eq:Q}
Q(u + v) = Q(u) + Q(v) + \Omega (u,v) \,\, (u,v \in V), \quad Q(b) = 1 \,\,
(b \in B) \ .
\end{equation}
(Clearly, these properties uniquely determine the restriction of $Q$ to $U$.)
An easy check shows that $Q(\tau_v (x)) = Q(x)$ whenever $Q(v) = 1$;
in particular, the function $Q$ is $\Gamma$-invariant.

Now everything is ready for a description of $\Gamma$-orbits in $V$. 

\begin{theorem}
\label{th:orbits}
If the graph $B$ is $E_6$-compatible then $\Gamma$ has precisely two orbits in
each set $(v + U) \setminus V^\Gamma$.
These two orbits are intersections of $(v + U) \setminus V^\Gamma$ with the
level sets $Q^{-1} (0)$ and $Q^{-1} (1)$ of $Q$.
\end{theorem}

The proof will be given in the next section. 
Let us show that this theorem implies Theorem~\ref{th:number of orbits} and 
Corollary~\ref{cor:number of components}.

\begin{corollary}
\label{cor:number of orbits}
If the graph $B$ is $E_6$-compatible then the number of $\Gamma$-orbits in 
$V$ is equal to 
$2^{\dim (V/U)} \cdot (2 + 2^{\dim (U \cap {\rm Ker} \ \Omega)})$;
in particular, if $U \cap {\rm Ker} \ \Omega = \{0\}$ then this number is
$3 \cdot 2^{\dim (V/U)}$.
\end{corollary}

\proof
By Proposition~\ref{pr:Gamma-invariants} and Theorem~\ref{th:orbits}, each
slice $v+ U$ with $\Omega (K,v + U) = 0$ splits into $2^{\dim K} + 2$ 
$\Gamma$-orbits, while each of the remaining slices splits into $2$ orbits.
There are $2^{\dim (V^\Gamma/K)}$ slices of the first kind and
$2^{\dim (V/U)} - 2^{\dim (V^\Gamma/K)}$ slices of the second kind.
Thus the number of $\Gamma$-orbits in $V$ is equal to
$$2^{\dim (V^\Gamma/K)} \cdot (2^{\dim K} + 2) + (2^{\dim (V/U)} - 2^{\dim
(V^\Gamma/K)}) \cdot 2 \ .$$
Our statement follows by simplifying this answer.
\endproof

Now Theorem~\ref{th:number of orbits} is just a reformulation of this 
Corollary. As for Corollary~\ref{cor:number of components}, one only needs 
to show that its assumptions imply that $U \cap {\rm Ker} \ \Omega = \{0\}$.
But this follows at once from Proposition~\ref{pr:boundary vertex}. 

\subsection{Proof of Theorem~\ref{th:orbits}}
We split the proof into several lemmas.
Let $E \subset U$ be the linear span of $6$ vectors from $B$ that form an
induced subgraph isomorphic to $E_6$.
The restriction of $\Omega$ to $E$ is nondegenerate; in particular, 
$E \cap K = \{0\}$.

\begin{lemma}
\label{lem:Q on E}
(a) Every $4$-dimensional vector subspace of $E$ contains at
least two non-zero vectors with $Q = 0$.

\noindent (b) Every $5$-dimensional vector subspace of $E$ contains at
least two vectors with $Q = 1$.
\end{lemma}

\proof
(a) It suffices to show that every $3$-dimensional subspace of $E$ contains
a non-zero vector with $Q = 0$.
Let $e_1, e_2$, and $e_3$ be three linearly independent vectors.
If we assume that $Q = 1$ on each of the $6$ vectors $e_1, e_2, e_3, e_1 +
e_2, e_1 + e_3$, and $e_2 + e_3$ then, in view of (\ref{eq:Q}), we must have
$\Omega(e_1, e_2) = \Omega(e_1, e_3) = \Omega(e_2, e_3) = 1$.
But then $Q(e_1 + e_2 + e_3) = 0$, as required.

(b) It follows from the results in~\cite{Janssen} (or by direct
counting) that $E$ consists of $28$ vectors with $Q = 0$ and $36$ vectors 
with $Q = 1$.
Since the cardinality of every $5$-dimensional subspace of $E$ is $32$, 
our claim follows.
\endproof

\begin{lemma}
\label{lem:Q nonconstant}
The function $Q$ is nonconstant on each set $(v + U) \setminus V^\Gamma$.
\end{lemma}

\proof
Suppose $v \in V \setminus V^\Gamma$.
By Lemma~\ref{lem:Q on E} (b), there exist two vectors $e \neq e'$ in $E$
such that
$$\Omega (v, e) = \Omega (v,e') = 0, \,\, Q(e) = Q(e') = 1 \ .$$
In view of (\ref{eq:Q}), we have $Q(v + e) = Q(v + e') = Q(v) + 1$, and it
is clear that at least one of the vectors $v + e$ and $v + e'$ is not 
$\Gamma$-invariant
(otherwise we would have $\Omega(e-e', u) = 0$ for all $u \in U$, which
contradicts the fact that $\Omega|_E$ is nondegenerate). 
\endproof

To prove Theorem~\ref{th:orbits}, it remains to show that $\Gamma$ acts
transitively on each level set of $Q$ in $(v + U) \setminus V^\Gamma$.
To do this, we shall need the following important result due to 
Janssen~\cite[Theorem~3.5]{Janssen}. 

\begin{lemma}
\label{lem:Q=1 transvections}
If $u$ is a vector in $U \setminus K$ such that $Q(u) = 1$ then the
symplectic transvection $\tau_u$ belongs to $\Gamma$.
\end{lemma}

We also need the following result from~\cite[Lemma~4.3]{ssv2}.

\begin{lemma}
\label{lem:Janssen}
If the graph $B$ is $E_6$-compatible then $\Gamma$ acts transitively on each of
the level sets of $Q$ in $U \setminus K$.
\end{lemma}

To continue the proof, let us introduce some terminology.
For a linear form $\xi \in U^*$, denote
$$T_\xi = \{u \in U \setminus K: Q(u) = \xi (u) = 1\} \ .$$
We shall call a family of vectors $(u_1, u_2, \dots, u_s)$ 
\emph{weakly orthogonal} if $\Omega(u_1 + \cdots + u_{i-1}, u_i) = 0$ 
for $i = 2, \cdots, s$.

\begin{lemma}
\label{lem:straight sums}
Let $\xi \in U^*$ be a linear form on $U$ such that $\xi|_K \neq 0$.
Then every nonzero vector $u \in U$ such that $Q(u) = \xi (u)$ can be
expressed as the sum $u = u_1 + \cdots + u_s$ of some weakly orthogonal 
family of vectors $(u_1, u_2, \dots, u_s)$ from $T_\xi$.
\end{lemma}

\proof
We need to construct a required weakly orthogonal family 
$(u_1, u_2, \dots, u_s)$ in each of the following three cases.

\noindent {\bf Case 1.} Let $0 \neq u = k \in K$ be such that 
$Q(k) = \xi(k) = 0$. Since $\xi \neq 0$, we have $\xi (b) = 1$ for some 
$b \in B$. By (\ref{eq:Q}), we also have $Q(b) = 1$.
Since $b \notin K$, we can take $(u_1, u_2) = (b, k-b)$ as a desired
weakly orthogonal family.

\noindent {\bf Case 2.} Let $u = k \in K$ be such that $Q(k) = \xi (k) = 1$.
By Lemma~\ref{lem:Q on E} (a), there exist distinct nonzero vectors $e$
and $e'$ in $E$ such that
$Q(e) = \xi (e) = Q(e') = \xi (e') = \Omega (e, e') = 0$.
Then we can take $(u_1, u_2, u_3) = (k - e, k - e', e + e' - k)$
as a desired weakly orthogonal family.

\noindent {\bf Case 3.} Let $u \in U \setminus K$ be such that 
$Q(u) = \xi(u) = 0$.
Since $\xi|_K \neq 0$, we can choose $k \in K$ so that $\xi (k) = 1$.
If $Q(k) = 1$ then a desired weakly orthogonal family for $u$ can be chosen as
$(u_1, u_2, u_3, u_4)$, where
$(u_1, u_2, u_3)$ is a weakly orthogonal family for $k$ constructed in 
Case 2 above, and $u_4 = u - k$.
If $Q(k) = 0$, choose $e \in E$ such that $Q(e) =1, \Omega (u,e) = 0$, and
$u - e \notin K$ (the existence of such a vector $e$ follows from
Lemma~\ref{lem:Q on E} (b)).
If $\xi (e) = 1$ then a desired weakly orthogonal family for $u$ can be 
chosen as $(u_1, u_2) = (e, u-e)$.
Finally, if $\xi (e) = 0$ then a desired weakly orthogonal family for 
$u$ can be chosen as $(u_1, u_2) = (e + k, u-e-k)$.
\endproof

Now everything is ready for completing the proof of Theorem~\ref{th:orbits}.
Take any slice $v + U \in V/U$; we need to show that $\Gamma$ acts
transitively on each of the level sets of $Q$ in $(v + U) \setminus V^\Gamma$.
First suppose that $(v + U) \cap V^\Gamma \neq \emptyset$; by
Proposition~\ref{pr:Gamma-invariants}, this means that $\Omega (K,v + U) = 0$.
Without loss of generality, we can assume that $v$ is $\Gamma$-invariant.
Then $\Omega (u, v) = 0$ for any $u \in U$, so we have 
$Q(v + u) = Q(v) + Q(u)$.
On the other hand, we have $g(v+u) = v + g(u)$ for any $g \in \Gamma$ and
$u \in U$.
Thus the correspondence $u \mapsto v+u$ is a $\Gamma$-equivariant bijection
between $U$ and $v + U$ preserving partitions into the level sets of $Q$.
Therefore our statement follows from Lemma~\ref{lem:Janssen}.

It remains to treat the case when $\Omega (K,v + U) \neq 0$.
In other words, if we choose any representative $v$ and define the linear
form $\xi \in U^*$ by $\xi (u) = \Omega (u, v)$ then $\xi|_K \neq 0$.
Let $u \in U$ be such that $Q(v) = Q(v + u)$; we need to show that $v+u$ 
belongs to the $\Gamma$-orbit $\Gamma (v)$.
In view of (\ref{eq:Q}), we have $Q(u) = \xi (u)$.
In view of Lemma~\ref{lem:straight sums}, it suffices to show that 
$\Gamma (v)$ contains $v + u_1 + \cdots + u_s$ for any weakly orthogonal 
family of vectors $(u_1, u_2, \dots, u_s)$ from $T_\xi$.
We proceed by induction on $s$.
The statement is true for $s =1$ because $v+ u_1 = \tau_{u_1} (v)$, and
$\tau_{u_1} \in \Gamma$ by Lemma~\ref{lem:Q=1 transvections}.
Now let $s \geq 2$, and assume that $v' = v + u_1 + \cdots + u_{s-1} \in
\Gamma (v)$.
The definition of a weakly orthogonal family implies that
$$v + u_1 + \cdots + u_s = v' + u_s = \tau_{u_s} (v') \in \Gamma (v) \ ,$$
and we are done.
This completes the proof of Theorem~\ref{th:orbits}.
\endproof

\section{Connected components of real double Bruhat cells}
\label{sec:components-claims}

In this section we give a (conjectural) geometric application of the above 
constructions. We assume that $\Pi$ is a Dynkin graph of simply-laced type, 
i.e., every connected component of $\Pi$ 
is the Dynkin graph of type $A_n, D_n, E_6, E_7$, or $E_8$.
Let $G$ be a simply connected semisimple algebraic group with 
the Dynkin graph $\Pi$. 
We fix a pair of opposite Borel subgroups $B_-$ and~$B$ in $G$;
thus $H=B_-\cap B$ is a maximal torus in~$G$. 
Let $N$ and $N_-$ be the unipotent radicals of $B$ and~$B_-$, respectively.
Let $\{\alpha_i: i \in \Pi\}$ be the system of simple roots for
which the corresponding root subgroups are contained in~$N$. 
For every $i \in \Pi$, let $\varphi_i: SL_2 \to G$ be 
the canonical embedding corresponding to $\alpha_i\,$. 
The (split) real part of $G$ is defined as the subgroup $G(\RR)$ of $G$ 
generated by all the subgroups $\varphi_i (SL_2(\RR))$. 
For any subset $L \subset G$ we define its real part by 
$L(\RR) = L \cap G(\RR)$. 

The \emph{Weyl group} $W$ of $G$ is defined by $W = {\rm Norm}_G (H)/H$. 
It is canonically identified with the Coxeter group $W(\Pi)$ (as defined in
Section~\ref{sec:Coxeter groups general}) via $s_i = \overline {s_i} H$, where
$$\overline {s_i} = \varphi_i \mat{0}{-1}{1}{0} \in {\rm Norm}_G (H) \ .$$ 

The representatives $\overline {s_i} \in G$ satisfy the braid relations in~$W$;
thus the representative $\overline w$ can be unambiguously defined for any 
$w \in W$ by requiring that 
$\overline {uv} = \overline {u} \cdot \overline {v}$
whenever $\l (uv) = \l (u) + \l (v)$.

The group $G$ has two \emph{Bruhat decompositions}, 
with respect to $B$ and $B_-\,$:
$$G = \bigcup_{u \in W} B u B = \bigcup_{v \in W} B_- v B_-  \ . $$
The \emph{double Bruhat cells}~$G^{u,v}$ are defined by 
$G^{u,v} = B u B  \cap B_- v B_- \,$.

Following \cite{BZ99}, we define the \emph{reduced double Bruhat cell}
$L^{u,v} \subset G^{u,v}$ as follows:
\begin{equation}
\label{eq: left sections}
L^{u,v} = N \overline u N  \cap B_- v B_- \ .
\end{equation}
The maximal torus $H$ acts freely on $G^{u,v}$ by left (or right) 
translations, and $L^{u,v}$ is a section of this action. 
Thus $G^{u,v}$ is biregularly isomorphic to $H \times L^{u,v}$, and all 
properties of $G^{u,v}$ can be translated in a straightforward way into 
the corresponding properties of $L^{u,v}$ (and vice versa). 
In particular, Theorem~1.1 in \cite{FZ} implies that $L^{u,v}$ is biregularly 
isomorphic to a Zariski open subset of an affine space of dimension 
$\l(u)+\l(v)$. 

\begin{conjecture}
\label{con:components}
{\rm For every two elements $u$ and $v$ in $W$, and every
reduced word $\ii \in R(u,v)$,
the connected components of $L^{u,v}(\RR)$ are in a natural bijection
with the $\Gamma_\ii (\FF_2)$-orbits in 
$\FF_2^{\l(u) + \l(v)}$ .}
\end{conjecture}

The precise form of this conjecture comes from the ``calculus of 
generalized minors" developed in \cite{FZ} and in a forthcoming 
paper \cite{BZ99}. 
If $u$ is the identity element $e \in W$ then $L^{e,v} = N \cap B_- v B_-$
is the variety $N^v$ studied in \cite{BZ}. 
When $G = SL_n$, and $v = w_0$, the longest element
in $W$, the real part $N^{w_0}(\RR)$ is the semi-algebraic set $N_n^0$ 
discussed in the introduction; in this case, the conjecture was proved  
in \cite{ssv1, ssv2} (for a special reduced word 
$\ii = (1, 2, 1, \dots, n-1, n-2, \dots, 2, 1) \in R(w_0)$).

\end{document}